\documentclass[11pt,reqno,oneside]{amsart}
\usepackage{latexsym}
 \usepackage{euscript}
 \renewcommand{\mathcal}{\CMcal}
\author{F.~Boniver, S.~Hansoul, P.~Mathonet, N.~Poncin}
\date{June 20, 2002}
\title[Equivariant Symbol Calculus for Forms]{Equivariant Symbol Calculus for Differential Operators 
Acting on Forms}
\newtheorem{lem}{Lemma}
\newtheorem{thm}[lem]{Theorem}
\newtheorem{prop}[lem]{Proposition}
\newtheorem{cor}[lem]{Corollary}

\theoremstyle{remark}

\theoremstyle{definition}

\newcommand{\R}{\mathbb{R}}
\newcommand{\PP}{\mathbb{P}}

\newcommand{\N}{\mathbb{N}}

\newcommand{\I}{\mathcal{I}}
\newcommand{\T}{\mathcal{T}}

\renewcommand{\L}{\mathcal L}
\newcommand{\D}{\mathcal{D}}
\renewcommand{\S}{\mathcal S}
\newcommand{\g}{\mathfrak{g}}

\newcommand{\euler}{\mathcal{E}}

\newcommand{\del}[1]{\frac{\partial}{\partial x^{{#1}}}}

\newcommand{\casinil}{N_\mathcal{C}}

\newcommand{\eval}[1]{\langle{#1}\rangle}
\newcommand{\restr}[1]{\mbox{}_{\vert {#1}}}
\newcommand{\alg}[2]{\mathit{{#1}}({#2})}

\newcommand{\casit}{C}
\newcommand{\casiop}{\mathcal{C}}
\newcommand{\algproj}{\mathit{sl}_{n+1}}
\newcommand{\symsp}{\mathcal{S}_p}
\newcommand{\im}{\mathrm{im}}

\begin{document}
    
    \begin{abstract}
        We prove the existence and uniqueness of a projectively equivariant 
        symbol map (in the sense of Lecomte and Ovsienko~\cite{LO}) for the spaces 
        $\D_p$ of differential operators transforming $p$-forms into 
        functions.  These results hold over a smooth manifold endowed with a flat projective structure.

As an application, we classify the $\mathit{Vect}(M)$-equivariant maps 
from $\D_p$ to $\D_q$ over any manifold $M$, recovering and 
improving earlier results by N.~Poncin~\cite{npo}.  This provides the complete answer to a question raised by P.~Lecomte about the extension of a certain intrinsic homotopy operator.
        \end{abstract}
        \maketitle
    \section{Introduction}
Let $\D^k_{\lambda}(M)$ denote the space of differential operators of order 
at most $k$ acting on $\lambda$-densities over a smooth manifold $M$. Let also $\D_{\lambda}(M)$ denote the filtered union of these spaces.
Both $\D^k_\lambda(M)$ and $\D_\lambda(M)$ are in a natural way modules over the Lie algebra $\mathit{Vect}(M)$ of vector fields 
over $M$.
If the manifold is furthermore endowed with a flat projective structure 
(i.e., can be locally identified with $\R\PP^n$, $(n=\dim M)$), the algebra 
$\mathit{sl}(n+1,\R)$, can be embedded in $\mathit{Vect}(M)$ as the algebra 
of infinitesimal projective transformations.  We denote the latter subalgebra by $\algproj$.

In \cite{LO}, P. Lecomte and V. Ovsienko considered the space $\D_{\lambda}(M)$ 
as a module over $\algproj$ and showed that it is actually 
isomorphic to the space of symbols $\mathit{Pol}(T^*M)$, which is the space 
of functions on $T^*M$ that are polynomial along the fibres.
Up to a natural normalization condition, the isomorphism $Q_{\lambda}$ from 
$\mathit{Pol}(T^*M)$ to $\D_{\lambda}(M)$ is unique and named the 
\emph{projectively equivariant quantization map}. Its inverse 
$\sigma_{\lambda}$ is the (projectively equivariant) \emph{symbol map}. An explicit 
formula is given in \cite{LO} for both mappings in terms of a divergence 
operator.  Furthermore, the knowledge of the equivariant symbol map proves useful to classify
 the $\mathit{Vect}(M)$-equivariant maps between $\D_\lambda(M)$ and $\D_{\mu}(M)$.

Several extensions of this work have recently been proposed.
In \cite{DLO}, one considers the spaces $\D_{\lambda,\mu}(M)$ of 
differential operators mapping $\lambda$-densities into 
$\mu$-densities as modules over the Lie algebra of infinitesimal 
conformal transformations, while in~\cite{bmIFFT}, the algebras under 
consideration correspond to the action of the symplectic (resp. 
pseudo-orthogonal) group on the Lagrangian (resp. pseudo-orthogonal) 
Grassmann manifolds. The results in both cases are the existence and 
uniqueness of an equivariant symbol map from $\D_{\lambda,\mu}(M)$ to 
the corresponding symbol space, provided that the value $\mu-\lambda$ 
does not belong to a countable and discrete critical set. Unfortunately, in these situations, a general formula 
for the symbol and quantization maps seems out of 
reach.

In this paper, we present a first example of projectively 
equivariant symbol calculus for differential operators acting on 
tensor fields.
We consider the spaces $\D_{p}$ of differential operators mapping $p$-forms into functions over a projectively flat manifold, whose dimension is assumed to be greater than 1. We still denote by $\mathcal{D}^k_p$ the subspace of operators with order at most $k$.
The corresponding symbol space $\mathcal{S}_{p}$ is made up of polynomial functions valued in contravariant 
antisymmetric tensor fields. 

We show that there exists a unique (up to normalisation) projectively 
equivariant symbol map from $\D_{p}$ to $\mathcal{S}_{p}$. We 
obtain an explicit formula in terms of the divergence operator and 
classical invariants of the space of symbols, such as the 
Koszul differential.
Next, we use this formula to classify the 
$\mathit{Vect}(M)$-invariant maps from $\D^k_{p}$ to $\D^l_{q}$ over \emph{any} 
manifold $M$ (not necessarily projectively flat anymore), such that $n=\dim M>1$.
We recover in this way the results by N. Poncin~\cite{npo} for $p\leq n-2$ 
and give the complete list of invariants for $p=n-1$ and $p=n$.

The paper is organized as follows.  In Section~\ref{sec.bas}, we recall basic definitions and notations.  In Section~\ref{sec.quant}, we collect some properties of the spaces of symbols and describe the projectively equivariant quantization map (see Theorem~\ref{exist}).  In Section~\ref{invop}, we determine the invariant operators from $\mathcal{S}_p$ to $\mathcal{S}_q$.  This is a first step towards the classification of $\mathit{Vect}(M)$-invariant maps from $\mathcal{D}^k_p$ to $\mathcal{D}^l_q$ over an arbitrary manifold $M$, which is completed in Section~\ref{sec.invopop}.      
\section{Basic definitions and notations}\label{sec.bas}
Throughout this paper, we assume manifolds to be second-countable, smooth, Hausdorff, and connected.
Let $M$ be such a manifold.
\subsection{Differential operators acting on $p$-forms}
Let $\Omega^p(M)$ denote the space of $p$-forms over $M$. 
The action of $\mathit{Vect}(M)$ on $\Omega^p(M)$ is standard:
\begin{equation}\label{LX}
    (L_{X}\omega)_{x}\,=\,(X.\,\omega)_x - \rho(D_{x}X)\,\omega_{x},
    \end{equation}
    where $X.$ denotes the differentiation along the vector 
    field $X$, $D_{x}X$ is the Jacobian matrix of $X$, and $\rho$ denotes the 
    natural action of $\mathit{gl}(n,\R)$ on the fibre $\wedge 
    ^pT^*_{x}M$.
    
    We denote by $\D^k_{p}$  the space of linear differential operators 
    of order at most $k$ from $\Omega^p(M)$ to 
    $\Omega^0(M)\,=\,C^{\infty}(M)$. The action of $\mathit{Vect}(M)$ 
    on this space is given by the commutator: if $D\,\in\,\D^k_{p}$, 
    \begin{equation}\label{LDX}
        \mathcal{L}_{X}D=L_{X}\circ D - D \circ L_{X}.
        \end{equation}
The filtered union $\D_p=\cup_k \D^k_p$ thus inherits a $\mathit{Vect}(M)$-module structure as well.
\subsection{Space of symbols}
For any $D$ in $\D^k_{p}(M)$, in local coordinates over a chart domain diffeomorphic to $\R^n$,
\begin{equation}\label{D}(D\,\omega)_{x}\,=\,\sum_{\mid\alpha\mid\leq k}\eval{A_{\alpha}(x),\,(\del{1})^{\alpha_{1}}\cdots 
(\del{n})^{\alpha_{n}}\,\omega},
\end{equation}
where $\eval{,}$ denotes the evaluation and  $A_{\alpha}\in\Gamma(\wedge^pT\R^n)$ for each $\alpha\in \N^n$.

The principal symbol $\sigma(D)$ of $D$ is the smooth section of the bundle 
\[
E^k_{p}\,=\,\wedge^pTM\otimes S^kTM\rightarrow M,
\]
 defined, for 
    all $x\in M$ and $\xi\in T^*_{x}M$ by 
\begin{equation}\label{sigma}
\sigma(D)_{x}(\xi)\,=\,\sum_{\mid\alpha\mid\ = k}\,(\xi_{1})^{\alpha_{1}}\cdots 
(\xi_{n})^{\alpha_{n}}\,A_{\alpha}(x).
\end{equation}
    We call the \emph{k-th order symbol space} and denote by 
    $\mathcal{S}^k_{p}$ the space $\Gamma(E^k_{p})$.
    The \emph{total symbol space} $\mathcal{S}_{p}$ is then the graded
    sum of these.
    Both $\mathcal{S}^k_{p}$ and $\mathcal{S}_{p}$ are modules over 
    the Lie algebra of vector fields. The Lie derivative is still given by~\eqref{LX}, $\rho$ now standing for  the natural action of 
    $\mathit{gl}(n,\R)$ on $\wedge^pT_{x}M\otimes S^kT_{x}M$.

\subsection{Projectively equivariant quantization and symbol maps}\label{subsec.quant}
Let $\g$ be a Lie subalgebra of $\mathit{Vect}(M)$. A \emph{$\g$-equivariant 
symbol map} is a $\g$-module isomorphism $\sigma_{\g}$: 
$\D_{p}(M)\rightarrow {\mathcal S}_{p}(M)$, such that for all 
$D\in\D^k_{p}(M)$, $$\sigma_{\g}(D) - \sigma(D)\in\bigoplus_{r=0}^{k-1}{\mathcal S}^{r}_{p}.$$
In analogy with [1,2], the inverse map $Q_{\g}$ will be called a \emph{$\g$-equivariant quantization 
map}.

Let us consider a basic example. We take $M=\R^n$ and consider the affine subalgebra 
$\mathit{Aff}$ of $\mathit{Vect}(\R^n)$, which is generated by constant and linear vector fields. The map
\begin{multline*}
\sigma_{\mathit{Aff}}: \D_{p}\rightarrow \S_{p}:\\
 \sum_{\mid\alpha\mid\leq k}\eval{A_{\alpha}(x),\,(\del{1})^{\alpha_{1}}\cdots 
(\del{n})^{\alpha_{n}}}\mapsto\sum_{\mid\alpha\mid \leq k}\,(\xi_{1})^{\alpha_{1}}\cdots 
(\xi_{n})^{\alpha_{n}}\,A_{\alpha}(x)
\end{multline*}
is a linear bijection. Moreover,
$$\sigma_{\mathit{Aff}}\circ \L_{X}\,=\,L_{X}\circ \sigma_{\mathit{Aff}}$$
for all $X\in\mathit{Aff}$.
In other words, the map $\sigma_{\mathit{Aff}}$ is an affinely 
equivariant symbol map.

Let now $M$ denote a manifold with a flat projective structure.  In the local coordinates defined by a suitable atlas, the infinitesimal projective transformations are generated by the vector fields
\begin{equation}
  \label{eq.sln1}
\del{r},\qquad x^s\,\del{r},\quad\mbox{and}\qquad x^r\euler,
\end{equation}
where $\euler$ denotes the Euler (or Liouville) vector field $\sum_{r=1}^nx^r\del{r}.$
These vector fields make up the \emph{projective Lie algebra}.  We denote it by $\algproj$ to recall that it is isomorphic to $\alg{sl}{n+1,\R}$.

It is well-known that $\algproj$ is a maximal subalgebra of the Lie algebra 
$\mathit{Vect}_{*}(\R^n)$ of vector fields with polynomial coefficients (see \cite{LO} for a proof).

Our approach to the existence and uniqueness problems for an $\algproj$-equivariant symbol map first involves some local computations.  Unless otherwise stated, we will from now on work in $\R^n$, with the local form~\eqref{eq.sln1} of $\algproj$.

\subsection{Polynomial formalism}
In order to simplify the computations involving differential operators, it is convenient to represent them by
polynomial functions belonging to the space of symbols, ``\`a la Fourier''.  We will still denote by $\L$ the Lie derivative of differential operators expressed in these symbolic terms.  
In other words, we equip $\S_p$ with a $\mathit{Vect}(\R^n)$-module structure isomorphic to that of $\D_p$, the isomorphism then being $\sigma_{\mathit{Aff}}$.

Seeking a $\g$-equivariant quantization map then amounts to seeking a $\g$-module 
isomorphism from $(\S_{p},L)$ to $(\S_{p},\L)$. 

Moreover, as in \cite{LMT}, we  consider the natural extension of 
$\sigma_{\mathit{Aff}}$ to multidifferential operators. We 
therefore represent such operators by polynomial functions of several 
variables.

{\bf Notation:} As we continue, we represent by $\eta$ the derivatives 
acting on the coefficients of an element of $\S_p$ and by $\zeta$ the derivatives acting on a vector 
field. We  also write the symbols as polynomial functions of the 
indeterminate $\xi$.

The formulas for the Lie derivatives can then be stated in a concise way.
\begin{lem} The Lie derivatives on the space of symbols $\S^k_{p}$ are given by 
the expressions
\[L_{X}u\,=\,\eval{X,\,\eta} u - \rho(X\otimes \zeta) u,\]
and 
\begin{equation}
  \label{eq.lieop}
\L_{X}u\,=\,L_{X}u + \sum_{r=1}^k t_{r}(X)\,u,
\end{equation}
where
\begin{multline}
  \label{eq.tr}
t_{r}(X): \S^k_{p}\rightarrow \S^{k-r}_{p}: \\ \Lambda\otimes 
P\mapsto - \Lambda\otimes X \frac{(\zeta 
\partial_{\xi})^{r+1}}{(r+1)!}P - X\wedge\,i_{\zeta}\Lambda\otimes \frac{(\zeta 
\partial_{\xi})^{r}}{r!}P.
\end{multline}
\end{lem}
\section{Projectively equivariant quantization map}\label{sec.quant}
Let us first review some properties of the space of symbols.
\subsection{Classical operators on the space of symbols}\label{subsec.class}
The typical fibre $W^k_{p}$ of $E^k_{p}$ is isomorphic 
    to $\wedge^p\R^n\otimes S^k\R^n$. 
The space $\bigoplus_{k,p} W^k_{p}$ is thus the tensor product of the exterior 
and symmetric algebras over $\R^n$.
Now, any $X\,\in\R^n$ 
can be identified with $X\otimes 1$ or $1\otimes X$. We
define two operators of multiplication by $X$:
\[\pi_{a}(X): W^k_{p}\rightarrow W^k_{p+1}: \Lambda\otimes P\mapsto 
X\wedge \Lambda\otimes P\]
and 
\[\pi_{s}(X): W^k_{p}\rightarrow W^{k+1}_{p}: \Lambda\otimes 
P\mapsto \Lambda\otimes XP.\]
We will also make use of the corresponding inner products by an 
element $\alpha\in\R^{n*}$:
\[i_{\alpha}: W^k_{p}\rightarrow W^k_{p-1}: \Lambda\otimes P\mapsto 
i_{\alpha}(\Lambda)\otimes P\]
and
\[\alpha\partial: W^k_{p}\rightarrow W^{k-1}_{p}: \Lambda\otimes 
P\mapsto \Lambda\otimes 
i_{\alpha}(P)\,=\,\Lambda\otimes\sum_{j=1}^n\alpha_{j}\partial_{\xi_{j}}P.\]
    Now, if $(v_{i})$ is a basis of $\R^n$ and if $(\beta^{i})$ is the dual 
    basis of $\R^{n*}$, the Koszul differential can be defined as
    \begin{equation}\label{delta}
        \delta: W^k_p\rightarrow W^{k-1}_{p+1}: u\mapsto 
        \sum_{j=1}^n\pi_{a}(v_{j})\circ (\beta^{j}\partial) (u),
        \end{equation}
        while the associated differential $\delta^*$ is defined by
        \begin{equation}\label{deltaet}
            \delta^*: W^k_p\rightarrow W^{k+1}{p-1}:u \mapsto 
        \sum_{j=1}^n i_{\beta^{j}}\circ\pi_{s}(v_{j})(u).
        \end{equation}
        The properties of these operators were discussed for instance in~\cite{Fe}. Let us quote that they commute with the natural action of 
        $\mathit{gl}(n,\R)$ on $W^k_p$, and that the following relations 
        hold on $W^k_p$:
\[
\delta^2=0,\quad (\delta^*)^2=0,\quad \delta\circ\delta^*+\delta^*\circ\delta=(k+p)\mathit{id}.
\]
            It follows that if $k+p\neq 0$, the mappings $\frac{\delta \circ
            \delta^*}{k+p}$ and $\frac{\delta^*\circ\delta}{k+p}$ are projectors 
            onto $\mathit{gl}(n,\R)$-submodules $A^k_p$ and $B^k_p$ of 
            $W^k_p$. At last, we know from the representation theory 
            of $\mathit{sl}(n,\R)$ that the decomposition 
            $$W^k_p\,=\,A^k_p \oplus B^k_p$$
            is the decomposition of $W^k_p$ into irreducible 
            $\mathit{sl}(n,\R)$-submodules.
It is worth noticing that, in terms of Young diagrams, $A^k_p$ is represented by a line of $k+1$
 boxes on top of $p-1$ lines with one single box, while $B^k_p$ is isomorphic to $A^{k-1}_{p+1}$.

Finally, the operators defined in this 
section can be canonically extended to the spaces $\S^k_{p}$. The operators 
$\delta$ and $\delta^*$ obtained in this way commute with the action of
$\mathit{Vect}(\R^n)$.  We write $\mathcal{A}^k_p=\im\,\delta\cap \mathcal{S}^k_p$ and $\mathcal{B}^k_p=\im\,\delta^*\cap \mathcal{S}^k_p$.
\subsection{Existence and uniqueness of the equivariant quantization}
In order to build the equivariant quantization, we use a method due to
Duval, Lecomte, and Ovsienko.  It appeared in~\cite{DLO}, in a conformal equivariance setting.    In~\cite{bmIFFT}, it proved to be relevant in a more general framework.  

This method relies on the comparison of eigenvalues and eigenvectors of two Casimir operators.  We denote by $C$ the Casimir operator associated to the representation $(\symsp,L)$ of $\algproj$ and by $\mathcal{C}$ the Casimir operator associated to the representation $(\symsp,\mathcal{L})$ of the same algebra.

Let us recall briefly the main stages of the process.

Assume that a $\algproj$-equivariant quantization $\mathcal{Q}$ is defined.  Then
\begin{equation}
  \label{eq:casi1}
\mathcal{Q}\circ\casit=\casiop\circ\mathcal{Q}.  
\end{equation}
 In particular, if $P$ is a homogeneous eigenvector of $C$ with eigenvalue $\alpha$, then $\mathcal{Q}(P)$ is an eigenvector of $\mathcal{C}$ with the same eigenvalue.  The principal symbol of $\mathcal{Q}(P)$ is $P$.

Conversely, assume that:
\begin{enumerate}
\item The Casimir operator $\casit$ is diagonalizable.
\item To any homogeneous eigenvector $P$ of $C$ with eigenvalue $\alpha$ corresponds a unique eigenvector of $\casiop$, which we denote $\mathcal{Q}(P)$, such that $\mathcal{Q}(P)$ has eigenvalue $\alpha$ and $P$ is the principal symbol of $\mathcal{Q}(P)$.
\end{enumerate}
Then the map $\mathcal{Q}$, linearly extended to $\symsp$, is the unique $\algproj$-equivariant quantization of $\symsp$.

The spectrum of $C$ is described by Proposition~\ref{pro.casival}.  The difference $\mathcal{C}-C$ is computed in Proposition~\ref{pro.casidiff}.  Finally, the unique $\algproj$-equivariant quantization is described in Theorem~\ref{exist}.
\begin{prop}\label{pro.casival}
The restriction of $\casit$ to $\mathcal{S}^k_p$ equals
\[
\frac{k+n+1}{n+1} \delta\circ\delta^* + \frac{k+n}{n+1}\delta^*\circ\delta.
\]
In particular, $C$ is diagonalizable and its spectrum is
\[
\{\overbrace{\frac{(k+n+1)(k+p)}{n+1}}^{\alpha^k_p}:k\in\N\}\cup\{\overbrace{\frac{(k+n)(k+p)}{n+1}}^{\beta^k_p}:k\in\N\}.
\]
\end{prop}
\begin{proof}
  We fix two bases of $\algproj$.  The first one is made up of the vector fields
\[
e_i=\del i, \euler, h_A=-\sum_{k,l} A^k_l x^l \del k, \epsilon^i=-\frac{1}{2(n+1)} x^i \euler,
\]
where $i\in\{1,\ldots,n\}$ and $A$ runs over a basis of $\alg{sl}{n,\R}$.
The second one is its dual with respect to the Killing form of $\algproj$.  We write
\[
\epsilon^i, \frac{1}{2n} \euler, \frac{n}{n+1} h_{A^*}, e_i,
\]
where $A^*$ is dual to $A$ with respect to the Killing form of $\alg{sl}{n,\R}$.

Then, 
\begin{equation}\label{eq.casidef}
C=2\sum_{i=1}^n L_{\epsilon^i}\circ L_{e_i}+ L_{\sum_{i} [e_i, \epsilon^i]}+ \frac{1}{2n} (L_\euler)^2+\frac{n}{n+1} \sum L_{h_A}\circ L_{h_{A^*}}.
\end{equation}
Since $C$ commutes with the constant vector fields, it has constant coefficients.  Hence, we  just  collect terms with such coefficients in the last right hand side.
From $L_X=X.-\rho(DX)$, where $\rho$ denotes the natural representation of $\alg{gl}{n,\R}$, we get
\[
C\restr{\mathcal{S}^k_p}=(\frac{k+p}{2}+\frac{(k+p)^2}{2n})\mathit{id}+\frac{n}{n+1}\sum_A \rho(A)\circ \rho(A^*).
\]
The second factor of the last term is the Casimir operator of $\alg{sl}{n,\R}$ acting on $W^k_p$.  It is a multiple of the identity on each irreducible component $A^k_p$ and $B^k_p$.  The computation of the eigenvalues is classical (see for instance~\cite[p.122]{hum}).

Since $\frac{1}{k+p}\delta\circ\delta^*$ and $\frac{1}{k+p}\delta^*\circ\delta$ are the projectors on $A^r_p$ and $B^r_p$, the conclusion follows.
\end{proof}
\begin{prop}\label{pro.casidiff}
  Let $\casinil$ denote the difference $\casiop-\casit$.  Then $\casinil$ equals
\[
\frac{1}{n+1}(\delta\circ(\eta\partial)\circ\delta^*+\delta^*\circ(\eta\partial)\circ\delta).
\]
\end{prop}
\begin{proof}
  We use the bases defined in the proof of Proposition~\ref{pro.casival}.  As a consequence of~\eqref{eq.lieop} and~\eqref{eq.casidef},
\[
\casinil=\mathcal{C}-C=2\sum_{i=1}^n t_1(\epsilon^i)\circ L_{e_i}.
\]
Using~\eqref{eq.tr}, we get
\[
\casinil\restr{\mathcal{S}^k_p}=\frac{k+p}{n+1}(\eta\partial)-\frac{1}{n+1}i_\eta\circ\delta.
\]
But, if $k>0$,
\[
\casinil\restr{\mathcal{S}^k_p}=\casinil\circ \frac{\delta\circ\delta^*+\delta^*\circ\delta}{k+p}.
\]

Since $[(\eta\partial),\delta]=0$ and $[(\eta\partial),\delta^*]=i_\eta$, the conclusion follows.
If $k=0$, $\casinil$ vanishes.
\end{proof}
\begin{thm}\label{exist}
The map $\mathcal{Q}$, defined by its restrictions
\[
\mathcal{Q}\restr{\mathcal{S}^k_p}=\mathrm{id}_{\mathcal{S}^k_p}+\sum_{l=1}^k \mathcal{Q}_{l,p},
\]  
with $\mathcal{Q}_{l,p}=\mathcal{Q}'_{l,p}+\mathcal{Q}''_{l,p}$, where
\[
\mathcal{Q}'_{l,p}=(\frac{1}{n+1})^l (\prod_{1\leq j\leq l}\frac{1}{\alpha^k_p-\alpha^{k-j}_p})(\delta\circ(\eta\partial)\circ\delta^*)^l,
\]
and
\[
\mathcal{Q}''_{l,p}=(\frac{1}{n+1})^l (\prod_{1\leq j\leq l}\frac{1}{\beta^k_p-\beta^{k-j}_p})(\delta^*\circ(\eta\partial)\circ\delta)^l,
\]
is the unique $\algproj$-equivariant quantization map of $\symsp$.
\end{thm}
\begin{proof}
  Let $P\in\symsp^k$ be an eigenvector of $C$.  The computation of $\mathcal{Q}(P)$ is similar according to whether $P\in\mathcal{A}^k_p$ or $P\in\im\,\mathcal{B}^k_p$ (cf. Subsection~\ref{subsec.class}).  Assume that $P\in \mathcal{A}^k_p$.  Then $C(P)=\alpha^k_p P$.  

The linear system
\begin{equation}
  \label{eq.casieq}
  \mathcal{C}(\mathcal{Q}(P))=\alpha^k_p \mathcal{Q}(P),
\end{equation}
which defines $\mathcal{Q}(P)$, is a triangular one.  Indeed, $C(\symsp^k)\subset \symsp^k$ and $\casinil(\symsp^k)\subset \symsp^{k-1}$.  Noticing that $\casinil$ stabilizes $\im\,\delta$, we write $\mathcal{Q}(P)=P+\sum_{r=1}^k P_{k-r}$, with $P_{k-r}\in\symsp^{k-r}$, and we observe that~\eqref{eq.casieq} amounts to
\[
P_{k-r}=\frac{1}{(n+1)(\alpha^k_p-\alpha^{k-r}_p)} \delta\circ(\eta\partial)\circ\delta^*(P_{k-r+1}).
\]
Hence the result.
\end{proof}
\section{Invariant differential operators from $\S^k_{p}$ to 
$\S^l_{q}$}\label{invop}
In this section, we will characterize the differential operators 
$\mathcal{T}$ from $\S^k_{p}$ to $\S^l_{q}$ that intertwine the action $L$
of $\algproj$ on these spaces. The arguments are 
essentialy those used in~\cite[Lemmas 10, 12]{bmIFFT}. We refer the reader to this 
paper for more detailed computations.
\begin{lem}
  If $\mathcal{T}: \S^k_{p}\rightarrow\S^l_{q}$ commutes with the 
  action of constant vector fields, then it has constant 
  coefficients. If in addition, it commutes with the action of the 
  Euler field, it is homogeneous of order $r\,=\,(k+p) -(l+q)$.
\end{lem}

\begin{thm}\label{thm.inv-s}
  The space of $\algproj$-invariant differential 
  operators from $\S^k_{p}$ to $\S^l_{q}$ is trivial if $(l,q)\not\in\{(k+1,p-1),(k,p),(k-1,p+1)\}$.
 The nontrivial spaces are generated by $\delta$,  $\mathit{id}$ and $\delta$, and $\delta$, respectively.
\end{thm}
\begin{proof}
  The assertion is trivial if $k+p=0$.  Assume that $k+p>0$ and let $\mathcal{T}$ denote an invariant mapping from $\S^k_{p}$ to $\S^l_{q}$.  We may 
  compose $\mathcal{T}$ on both sides with the projectors on $\im \, \delta$ and 
  $\im\, \delta^*$. 
  We obtain in this way (at most) four invariant operators. Let us 
  show for instance how one can characterize an invariant operator from
  $\mathcal{A}^k_{p}$ to 
  $\mathcal{B}^l_{q}$. 
  
  Since $\mathcal{T}$ has constant coefficients 
  and is homogeneous with order $r\,=\,(k+p) -(l+q)$, its 
  restriction to the subspace of $\mathcal{A}^k_{p}$ of sections with degree $r$ polynomial coefficients has values in the subspace of constant sections of $\mathcal{B}^l_{q}$.
  As modules over the Lie algebra of linear and divergence free vector 
  fields, these subspaces 
  identify with $S^r\R^{n*}\otimes A^k_{p}$ and $B^l_{q}$ 
  respectively, endowed with the standard representation of $\mathit{sl}(n,\R)$.
  
  If $\mathcal{T}$ is not trivial, since $B^l_{q}$ is irreducible, the decomposition of $S^r\R^{n*}\otimes A^k_{p}$ into 
  $\mathit{sl}(n,\R)$-irreducible submodules contains a factor isomorphic to $B^l_q$. We can compute the Young diagrams 
  associated to these submodules using Littlewood-Richardson rule and 
  compare these to the diagram associated to $B^l_{q}$. We then observe 
  that the diagram associated to $B^l_{q}$ is obtained by deleting $r$
  boxes in the diagram associated to $A^k_{p}$.
  This amounts to the information $l\leq k+1$ and $q+1\leq p$. 
  
  Since $\mathcal{T}$ commutes with the Casimir operator $C$, we also have 
  $\alpha^k_{p}=\beta^l_{q}$, and thus
\begin{eqnarray*}
    0=(n+1)(\alpha^k_{p} - \beta^l_{q})  &=&  
    (k+p)(k+n+1)-(l+q)(l+n)
\\&\geq &(l+q)(k-l+1)\geq 0.
  \end{eqnarray*}
This implies successively $l+q>0$, $l=k+1$, $p=q+1$, and finally $r=0$. Therefore $T$ defines a $\mathit{sl}(n,\R)$-invariant map 
  from $A^k_{p}$ to $B^{k+1}_{p-1}$. The space of such invariants is by 
  Schur's lemma at most one dimensional, and is thus generated by the 
  map $\delta^*$.
\end{proof}
\section{Invariant operators from $\D_{p}^k$ to $\D_{q}$}\label{sec.invopop}
In this section, we let $M$ denote an arbitrary manifold (not necessarily projectively flat). We show how the projectively equivariant quantization map allows us to classify the $\mathit{Vect}(M)$-invariant operators $T$ from
$\D_{p}^k(M)$ to $\D_{q}(M)$, i.e. those $T$ such that
\begin{equation}\label{eqinv}
  T(\L_{X}D)\,=\,\L_{X}T(D),
\end{equation}
for all $X\in\mathit{Vect}(M)$ and $D\in \D_{p}^k(M)$.

We denote by $\I^{k}_{p,q}$ the space of these invariant 
operators.

The classification method presented here allows us to recover quite easily the results of~\cite{npo} in the case 
$p\leq n-2$ ($n=\dim\,M$) and to complete the 
classification of invariant maps in the cases $p=n-1$ and $p=n$.
    
    The following result is proved in~\cite{npo}.
    \begin{prop}
        If $T\in\I^{k}_{p,q}$, then $T$ is a local operator.
        \end{prop}
    We can thus consider the restrictions of invariant operators to 
    relatively compact chart domains over the manifold. We can 
    furthermore consider chart domains that are diffeomorphic to 
    $\R^n$, and it will be sufficient to compute the spaces 
    $\I^{k}_{p,q}$ over the manifold $M=\R^n$.

    Now, we use the projectively equivariant quantization map as a tool in the following way.
 To any $T\in\I^{k}_{p,q}$, we  associate the operator
    $$\T: \bigoplus_{r=0}^k\S^r_{p}\rightarrow 
    \S_{q}: u\mapsto (\mathcal{Q}^{-1}\circ T\circ 
    \mathcal{Q})(u).$$
   We also denote by $\T_{r,p}$ the restriction of $\T$ to $\S^r_{p}$.
   In this framework, equation (\ref{eqinv}) means that $\T$ has to 
   commute with the operator 
   $\mathcal{Q}^{-1}\circ\L_{X}\circ\mathcal{Q}$, for all 
   $X\in\mathit{Vect}(\R^n)$. Restricted to $\S^r_p$, the latter operator 
 equals
   \[ L_{X} + 
   \sum_{i=1}^r \gamma_{i,p}(X),\]
   where each $\gamma_{i,p}(X): \S^r_{p}\rightarrow\S^{r-i}_{p}$ is a 
   differential operator with respect to $X$ and vanishes for all $X$ in 
   $\algproj$.
   Equation (\ref{eqinv}) then first implies 
   $$L_{X}\circ \T_{r,p} = \T_{r,p}\circ L_{X},$$ 
   for any $X\in\algproj$, and Theorem~\ref{thm.inv-s} ensures that no more than three nontrivial invariant types exist:
   \begin{itemize}
       \item If $q=p-1$, the map $\T_{r,p}: \S^r_{p}\rightarrow 
       \S^{r+1}_{p-1}$ writes $a_{r,p}\delta^*$ for some constant 
       $a_{r,p}.$
       \item If $q=p$, the map $\T_{r,p}: \S^r_{p}\rightarrow 
       \S^{r}_{p}$ writes $b_{r,p}\delta\circ\delta^* + 
       c_{r,p}\delta^*\circ\delta$ for some constants
       $b_{r,p}$ and $c_{r,p}$.
       \item If $q=p+1$, the map $\T_{r,p}: \S^r_{p}\rightarrow 
       \S^{r-1}_{p+1}$ writes $d_{r,p}\delta$ for some constant 
       $d_{r,p}.$
       \end{itemize}
       
       The commutation relation (\ref{eqinv}) then forces $\mathcal{T}$ to fulfil
       \begin{equation}\label{eqgamma}
         \gamma_{i,q}(X)\circ 
         \T_{r,p}\,=\,\T_{r-i,p}\circ\gamma_{i,p}(X),
       \end{equation}
       for all $r\leq k$, $i\leq r$, and $X\in\mathit{Vect}(\R^n)$.

An important remark helps simplifying the computations further at this stage.  Notice that it is sufficient to require that condition~\eqref{eqinv}, and therefore~\eqref{eqgamma}, be satisfied for all $X$ with polynomial coefficients.  Indeed, the Lie derivatives are differential operators with respect to $X$.  But then, since $\algproj$ is a maximal subalgebra of $\mathit{Vect}_*(\R^n)$ (cf.~Subsection~\ref{subsec.quant}), it suffices to verify that~\eqref{eqgamma} holds for a given polynomial vector field $X$ out of $\algproj$ in order to ensure that $T$ be a $\mathit{Vect}(\R^n)$-invariant.

\noindent\textbf{Notation:} From now on, we  therefore fix $X=(x^1)^2\del 2$.
       \begin{lem}\label{lem.comm}
           The following relations hold on $\S^r_{p}$:
           \begin{itemize}
               \item For all $i$, $r$, and $p$: $\delta^*\circ \gamma_{i,p}(X) = 
               \gamma_{i,p-1}(X)\circ\delta^*$.
               \item For all $r\geq 1$ and $1\leq p\leq n-1$, $\delta^*\circ 
               \gamma_{1,p}(X) \neq 0.$
               \item The operator $\gamma_{1,p}$ vanishes iff $r=0$ or $p=n$ or $(r,p)=(1,0)$.
               \item For all $r\geq 3$ and $p\in\{0,\ldots,n-2\}$,
               $\delta\circ\gamma_{1,p}(X)\circ\delta\neq 0$.
               \item For all $r\geq 2$, $\delta^*\circ\gamma_{2,n}(X)\neq 0$.
           \end{itemize}
           \end{lem}
           \begin{proof}
           The proof is straightforward. Let us perform it for the first claim. 
           Note that the maps $\gamma_{i,p}$ are 
           polynomial functions of the operators $L_{X}$, $t_{1}(X)$ and $Q_{j,p}$ 
           $(1\leq j\leq i)$. For instance, we have 
           \[\gamma_{1,p}(X)\,=\,t_{1}(X) + [L_{X},Q_{1,p}].\]
           In order to conclude, we notice that 
           $[\delta^*,L_{X}]\,=\,[\delta^*,t_{1}(X)]=0$ and that 
           $\delta^*\circ Q_{i,p}\,=\,Q_{i,p-1}\circ \delta^*$.
           \end{proof}     
           \begin{cor}\label{iso} For all $p$,   
           $\D_{p}^0$ is isomorphic to $\S^0_{p}$ as a $\mathit{Vect(M)}$-module, while
           $\D^1_{p}$ is isomorphic to $\S^1_{p}\oplus\S^0_{p}$, for $p=0$ and 
           $p=n$.
               \end{cor}
               \begin{proof} Over any chart domain $U$ of $M$, Theorem 
               \ref{exist} allows to define a unique isomorphism $\mathcal{Q}_{U}$  of 
               $\algproj$-modules from the 
               space of symbols to the space of differential operators over 
               the chart. 
               By lemma~\ref{lem.comm}, the operator $\gamma_{1,p}(X)$ vanishes and the maps 
               $\mathcal{Q}_{U}$ actually commute with the action of all 
               vector fields over the chart. 
               The uniqueness of the maps $\mathcal{Q}_{U}$ as $\mathit{Vect}(U)$-module 
               isomorphisms then implies that they glue toghether to define a 
               global isomorphism.
                   \end{proof}
           Now, before turning to the classification results, let us recall 
           the expression of some classical invariant maps and introduce some 
           new ones.
           \subsection{Some invariant maps}
           The dual $d^*$ of the de Rham differential is defined by
           $$d^*:\D^k_{p}\rightarrow\D^{k+1}_{p-1}: D\mapsto D\circ d.$$
           It is well-known that it commutes with the action of 
           $\mathit{Vect}(M)$.
           
           The space $\D_{0}$ is the direct sum of the 
           space of differential operators vanishing on the constants and 
           of $\D^0_{0}$. We denote by $\I_{0}$ the projection onto the 
           second summand:
           \[I_{0}: \D_{0}\rightarrow\D^0_{0}: D\mapsto D(1),\]
           which is of course invariant.

           The \emph{conjugation operator} was for instance presented in~\cite[5.5.3]{Math} in the framework of differential operators acting on densities.  
           Let $U$ be a chart domain of $M$.  Densities with weight $0$ over $U$ are nothing but functions, while $1$-densities identify with $n$-forms.  Here, the conjugation operator 
\[
\EuScript{C}: \mathcal{D}_n^k\to \mathcal{D}_n^k
\]
is  locally defined by the condition
\[
\int_U D(f)g=\int_U f \EuScript{C}(D)(g),
\]
where $f$ and $g$ are arbitrary compactly supported $n$-forms over $U$ and $D\in\mathcal{D}^k_n$ has also compact support in $U$.
It turns out that $\EuScript{C}$ is then a globally defined invariant.

    Corollary \ref{iso} allows to define a new invariant map, acting 
    on $\D^k_{p}$ when $k=1$ and $p\geq 1$ or $(k,p)=(2,n-1)$: 
    \[K: \D^k_{p}\rightarrow \D^{k-1}_{p+1}: D\mapsto 
    \mathcal{Q}\circ \delta\circ \sigma (D).\]
   
    Finally, we recall the existence of an invariant map $K':\D^2_{0}\rightarrow 
    \D^1_{1}$, which is defined by its local expression over a chart domain of $M$:
  \[ 
  \begin{array}{lll}\mathcal{Q}^{-1}\circ 
    K'\circ\mathcal{Q}\restr{\S^2_{0}}\,=\frac{1}{2}\delta, & \mathcal{Q}^{-1}\circ 
    K'\circ\mathcal{Q}\restr{\S^1_{0}}\,=\delta, &K'\restr{\D^0_{0}}=0.
   \end{array}
   \]
   P. Lecomte presented this operator in a more general setting as a homotopy operator for $d^*$ and proved its global existence and invariance (see \cite[p.~188]{leglob}).  
\subsection{Classification}
We now state the classification of all invariant maps from $\mathcal{D}^k_p$ to $\mathcal{D}_q$.  Theorem~\ref{thm.inv-op-ppm1} takes into account the case $q=p-1$.  Theorem~\ref{thm.inv-op-pp1} and Theorem~\ref{thm.inv-op-pp2} give the results for $q=p$. Only a few values of $k$ and $p$ yield invariants when $q=p+1$.  They are given, as well as the corresponding invariants themselves, by Theorem~\ref{thm.inv-op-ppp1}.
\begin{thm}\label{thm.inv-op-ppm1}
  If $p<n$ or $k=0$, then $\mathcal{I}_{p,p-1}^k$ is generated by $d^*$.  If $p=n$ and $k\neq 0$, then $\dim(\mathcal{I}_{p,p-1}^k)=2$, an additional generator being given by $d^*\circ\EuScript{C}$.
\end{thm}
\begin{thm}\label{thm.inv-op-pp1}
  Let $1\leq p\leq n-1$.  If $k\neq 1$ and $(k,p)\neq (2,n-1)$, then $\mathcal{I}_{p,p}^k$ is generated by the identity map.  If $k=1$ or $(k,p)=(2,n-1)$, then another generator is given by $d^*\circ K$.  
\end{thm}
\begin{thm}\label{thm.inv-op-pp2}
  The space $\mathcal{I}_{0,0}^0$  (resp. $\mathcal{I}_{n,n}^0$) is generated by $\mathit{id}$.  If $k\neq 0$, $\mathcal{I}_{0,0}^k$ (resp. $\mathcal{I}_{n,n}^k$) is generated by $\mathit{id}$ and $I_0$ (resp. $\EuScript{C}$).
\end{thm}
\begin{thm}\label{thm.inv-op-ppp1}
  If $k\neq 1$ and $(k,p)\not\in\{(2,0),(2,n-1)\}$, then $\mathcal{I}_{p,p+1}^k=\{0\}$.  If $k=1$ or $(k,p)=(2,n-1)$, the invariants are multiples of $K$.  If $(k,p)=(2,0)$, they are multiples of $K'$.
\end{thm}
\begin{proof}[Proof of Theorem~\ref{thm.inv-op-ppm1}]
Let $T\in \mathcal{I}_{p,p-1}^k$.  We follow the method outlined and use the notations introduced at the beginning of this Section.  We have $\mathcal{T}_{r,p}= a_{r,p} \delta^*$, for all $r\in\{0,\ldots,k\}$.  The conclusion follows if $k=0$.

Otherwise, the first of the commutation relations~\eqref{eqgamma} then yields
\[
a_{r,p} \gamma_{1,p-1}(X)\circ \delta^*\restr{\mathcal{S}^r_p} = a_{r-1,p} \delta^*\circ\gamma_{1,p}(X)\restr{\mathcal{S}^r_p}. 
\]
In view of Lemma~\ref{lem.comm}, we get
\[
(a_{r,p}-a_{r-1,p}) \delta^*\circ \gamma_{1,p}(X)\restr{\mathcal{S}^r_p}=0.
\]
If $p\neq n$, $\delta^*\circ \gamma_{1,p}(X)\restr{\mathcal{S}^r_p}\neq 0$ and $\dim(\mathcal{I}_{p,p-1}^k)\leq 1$.  Since $d^*$ is a suitable invariant, the conclusion follows.

If $p=n$, the second relation of~\eqref{eqgamma} yields
\[
(a_{r,p}-a_{r-2,p}) \delta^*\circ\gamma_{2,n}(X)\restr{\mathcal{S}^r_n}=0,
\]
for all $r\geq 2$. Lemma~\ref{lem.comm} then ensures that $\dim(\mathcal{I}^k_{n,n-1}) = 2$.  Hence the result. 
\end{proof}
\begin{proof}[Proof of Theorem~\ref{thm.inv-op-pp1}]
  Let $T\in \mathcal{I}^k_{p,p}$.  Then $d^*\circ T\in \mathcal{I}^{k}_{p,p-1}$ is an invariant operator classified by Theorem~\ref{thm.inv-op-ppm1} and thus a multiple of $d^*$.  If we add a suitable multiple of $\mathit{id}$ to $T$, we may assume that 
\[
 d^*\circ T=0. 
\]

Writing $\mathcal{T}_{r,p}=b_{r,p} \delta\circ \delta^*+c_{r,p} \delta^*\circ\delta$, we get $(r+p)b_{r,p}\delta^*\restr{\mathcal{S}^r_p}=0$ and thus $b_{r,p}=0$, since $p>0$.  The first relation of~\eqref{eqgamma} yields
\[
c_{r,p} \gamma_{1,p}(X)\circ\delta^*\circ\delta\restr{S^r_p}=c_{r-1,p} \delta^*\circ\delta \circ \gamma_{1,p}(X)\restr{\mathcal{S}^r_p}
\]
 and then
\begin{equation}\label{eq.dstar2}
\delta^*\circ(c_{r,p} \gamma_{1,p+1}(X)\circ\delta - c_{r-1,p} \delta\circ\gamma_{1,p}(X))\restr{\mathcal{S}^r_p}=0,
\end{equation}
for all $r\in\{1,\ldots,k\}$.
If $k=1$, this equation is trivial and $\dim(\mathcal{I}^1_{p,p})\leq 2$.  If $k>1$, we evaluate the left hand side on $\im\, \delta$ and obtain 
\[
c_{r-1,p} \delta^*\circ(\delta\circ\gamma_{1,p}(X)\circ\delta)\restr{\mathcal{S}^{r+1}_{p-1}}=0.
\] 
Since $\im\, \delta \cap\ker\,\delta^*=\{0\}$, we conclude that $\mathcal{T}=\mathcal{T}_{k,p}=c_{k,p} \delta^*\circ\delta$.  Another look at  $\eqref{eq.dstar2}$ shows that $c_{k,p}=0$ if $p\neq n-1$.  But, if $p=n-1$ and $k>2$, the second relation of~\eqref{eqgamma} forces $c_{k,p}=0$.  Hence the result.
\end{proof}
\begin{proof}[Proof of Theorem~\ref{thm.inv-op-pp2}]
We proceed as in the proof of the last two theorems.

We know that $\mathcal{T}_{r,p}=e_{r,p} \mathit{id}$, for $p\in\{0,n\}$.
The first (resp. second) relation of~\eqref{eqgamma} yields the result when $p=0$ (resp. $p=n$).
\end{proof}
\begin{proof}[Proof of Theorem~\ref{thm.inv-op-ppp1}]
Let $T\in \mathcal{I}^k_{p,p+1}$.
  We have $\mathcal{T}_{r,p}=d_{r,p} \delta$ and thus $\mathcal{T}_{r,p}=0$.  The result for $k=1$ is then clear. 

From Theorem~\ref{thm.inv-op-pp1}, we know that $d^*\circ T$ is a multiple of the identity, provided that $1\leq p\leq n-2$.  This is also true if  $k\geq 3$ and $p=n-1$. But the restriction of $d^*\circ T$ to operators of order $0$ vanishes.  Therefore, 
\[
d^*\circ T=0.
\]
We deduce that $d_{r,p} \delta^*\circ\delta\restr{\mathcal{S}^r_p}=0$ and thus $T=0$.

Now, if $p=n-1$ and $k=2$, we get similarly that $d^*\circ T$ is a multiple of $d^*\circ K$.

If $p=0$ and $k\geq 2$, $d^*\circ T$ is a combination of the identity and $I_0$.  It vanishes on operators of order $0$.  Therefore, $d^*\circ T=a(\mathit{id}-I_0)$ for some $a\in\R$. This equality implies $d_{r,0} \delta^*\circ\delta=a\,\mathit{id}$ on $\mathcal{S}^r_0$, for all $r\geq 1$, which in turn proves the result for $k=2$.
If $k\geq 3$, the first relation of~\eqref{eqgamma} forces once more $T$ to vanish.  Hence the result.
\end{proof}
   
\section*{Acknowledgments}
We would like to thank P.~Lecomte for valuable suggestions about this work.
The first two authors thank the Belgian National Fund for Scientific Research for their Fellowships.  The last author thanks the Luxembourgian Research Department for Grant MEN/C.U.L./99/007.
    \bibliographystyle{plain}
    \bibliography{Articleform}
Fabien Boniver, Sarah Hansoul, Pierre Mathonet, Universit\'e de Li\`ege,\\ D\'epartement de Math\'ematique, B. 37, Grande Traverse, 12, B-4000 Li\`ege, Belgium\\
e-mail: f.boniver@ulg.ac.be, s.hansoul@ulg.ac.be, p.mathonet@ulg.ac.be. \\[1em]
Norbert Poncin, Centre Universitaire de Luxembourg, D\'epartement de \\ Math\'ematiques, 162A, Avenue de la Fa\"\i encerie, L-1511 Luxembourg,\\ Luxembourg.  \\
e-mail: poncin@cu.lu
 \end{document}